\documentclass[a4paper, 10pt]{article}

\usepackage{amsmath, amssymb, amsthm, amsfonts}
\usepackage[utf8]{inputenc}
\usepackage[T2A, T1]{fontenc}
\usepackage{lmodern}
\usepackage[german, french, russian, UKenglish]{babel}

\usepackage{bookmark}
\hypersetup{colorlinks=true, linkcolor=blue, citecolor=blue}

\usepackage{csquotes}
\usepackage[backend=biber,%defernumbers=true,%
giveninits=false,isbn=false,doi=false,url=false,%
style=alphabetic,sorting=nyt,%
backref=false,
labelalpha=true,
bibencoding=utf8]{biblatex}
\usepackage[inline, shortlabels]{enumitem}
\usepackage[a4paper]{geometry}
\usepackage{comment}

\newcommand{\bF}{\mathbb{F}}
\newcommand{\bK}{\mathbb{K}}
\newcommand{\bN}{\mathbb{N}}
\newcommand{\bZ}{\mathbb{Z}}
\newcommand{\fp}{\mathfrak{p}}
\newcommand{\Ann}{\operatorname{Ann}}
\newcommand{\Aut}{\operatorname{Aut}}
\newcommand{\End}{\operatorname{End}}
\newcommand{\DefAut}{\operatorname{DefAut}}
\newcommand{\DefEnd}{\operatorname{DefEnd}}
\newcommand{\Fr}{\operatorname{Fr}}
\newcommand{\Frac}{\operatorname{Frac}}
\newcommand{\GL}{\operatorname{GL}}
\newcommand{\im}{\operatorname{im}}
\newcommand{\Id}{\operatorname{Id}}
\newcommand{\op}{\mathrm{op}}

\newcommand{\Stab}{\operatorname{Stab}}
\newcommand{\Vect}[1]{#1\text{-}\mathbf{Vect}}

\newcommand{\generated}[1]{\left\langle#1\right\rangle}

\theoremstyle{plain}
\newtheorem*{theorem*}{Theorem}
\newtheorem{corollary}{Corollary}
\newtheorem*{lemma*}{Lemma}
\newtheorem*{definition*}{Definition}
\newtheorem*{notation}{Notation}
\newtheorem{proposition}{Proposition}

\theoremstyle{definition}

\theoremstyle{remark}

\title{Zilber's skew-field lemma}
\author{Adrien Deloro}

\AtBeginBibliography{\small}
\renewbibmacro{in:}{%
\ifentrytype{article}{}{\printtext{\bibstring{in}\intitlepunct}}}

\begin{document}

\maketitle

\begin{center}
\begin{minipage}{.8\textwidth}
%\abstract{\noindent
\small
\textbf{Abstract.}
We revisit one of Zilber's early results in model-theoretic algebra, viz.~definability in Schur's lemma. This takes place in a broader context than the original version from the seventies.
\emph{The present exposition contains results extracted from more general research in progress with Frank O.~Wagner.}
%}
\end{minipage}
\end{center}
\normalsize
% \vspace{1cm}

\centerline{\rule{0.7\linewidth}{.5pt}}\medskip

\hfill
\begin{minipage}[t]{0.4\textwidth}
\small
\begin{flushright}
{\centering\itshape
La droite laisse couler du sable.\\
Toutes les transformations sont possibles.\\}
Paul Éluard
\end{flushright}
\end{minipage}
\bigskip

The present contribution discusses and proves a linearisation result originating in Zilber's early work.
(\begin{enumerate*}
\item 
$o$-minimal dimension and Borovik-Morley-Poizat rank are examples of finite dimensions.
\item
I have prefered not to conflate $T$ with $\bK$ in the statement.
\item
There are classical corollaries in \S~\ref{s:corollaries}.
\item
The result bears no relationship to indecomposable generation discussed in \S~\ref{s:indecomposable}.
\end{enumerate*})

\begin{theorem*}[Zilber's skew-field lemma]
Work in a finite-dimensional theory.
Let $V$ be a definable, connected, abelian group and $S, T \leq \DefEnd(V)$ be two invariant rings of definable endomorphisms such that:
\begin{itemize}
\item
$V$ is irreducible as an $S$-module;
\item
$C(S) = T$ and $C(T) = S$, with centralisers taken in $\DefEnd(V)$;
\item
$S$ and $T$ are infinite;
\item
$S$ \emph{or} $T$ is unbounded.
\end{itemize}
Then there is a definable skew-field $\bK$ such that $V \in \Vect{\bK}_{<\aleph_0}$; moreover $S \simeq \End(V\colon \Vect{\bK})$ and $T \simeq \bK \Id_V$ are definable.
\end{theorem*}

% \medskip
\centerline{\rule{0.7\linewidth}{.5pt}}\medskip

\noindent\S~\ref{S:introduction} %contains (amateurish and fortunately optional) historical remarks.
provides context.
\S~\ref{S:theorem} discusses the statement, and gives all definitions. The proof is in \S~\ref{S:proof}.

\section{Introduction}\label{S:introduction}

\S~\ref{s:Schur} explains the relation to Schur's Lemma. \S~\ref{s:history:short} makes some historical remarks. \S~\ref{s:soluble} discusses a more famous corollary on fields in abstract groups.

\subsection{Schur's Lemma}\label{s:Schur}

Among the early work of Zilber's are a couple of gems in model-theoretic algebra. The present article deals with one of the phenomena he discovered: \emph{many $\aleph_1$-categorical groups interpret infinite fields}.
The result, or the method, or the general line of thought, is often called \emph{Zilber's field theorem}.
It 
%is a definable analogue of the Artin-Wedderburn theorem, with lineage stemming
stems
 from Schur's lemma in representation theory.

\begin{lemma*}[Schur's Lemma]
Let $R$ be a ring and $V$ be a simple $R$-module. Then the covariance ring $\bF = C_{\End(V)}(R)$ is a skew-field, $V$ is a vector space over $\bF$, and $R \hookrightarrow \End_R(V)$.
\end{lemma*}

Zilber's deep observation is simple:
\begin{quote}
\itshape
in many model-theoretically relevant cases, $\bF$ is definable.
\end{quote}

A precise and modern form of the latter statement, given as Corollary~\ref{c:SZ}, is a straightforward consequence of the main theorem above.
(One should remember that every module is actually a bimodule by introducing Schur's covariance ring.)
I shall henceforth call it (in long form) the \emph{Schur-Zilber skew-field lemma}, hoping that Boris will not mind being in good company.
Far be it from me to minimise its significance by dubbing it a lemma instead of a theorem; quite the opposite as lemmas are versatile devices---methods.

\subsection{Editorial fortune of the lemma}\label{s:history:short}%, short}

This subsection is a layman's attempt at providing historical remarks. I apologise for misconceptions.

\begin{itemize}
\item
As one learns from \cite[p.~139]{CPioneers}, Schur's Lemma itself appears in \cite[\S~2, I.]{SNeue} with comment: `\emph{der auch in der \textsc{Burnside}'schen Darstellung der Theorie eine wichtige Rolle  spielt}'.
\item
Before Zilber's result was known, Cherlin \cite[4.2, Theorem~1]{CGroups} found a definable field independently. There interpretation is obtained by hand (and seemingly by miracle), without a general method. Cherlin heard about Zilber's work after completing his own; \cite[{}1.4]{CGroups} is very informative.
\item
The lemma itself seems not to have drawn as much attention as its corollary on soluble groups (\S~\ref{s:soluble}). There are few traces of the lemma as a stand-alone statement.
%TODO so... where?
\item
All sources discussing the topic \cite{ZGruppi, TClassification, ZSome, NSolvable, PGroupes, NNonassociative, LWCanada, BNGroups, MPPrimitive} rely on indecomposable generation (however see~\S~\ref{s:indecomposable}).
\item
This is different in the $o$-minimal context, but \cite[Theorem~2.6]{PPSSimple} has its own techniques. (The earlier \cite[Proposition~2.4]{NPRGroups}, which bears no reference to Zilber, resembles Cherlin's coordinatisation by hand \cite{CGroups}.)

This and the above ad may have given the impression that the Schur-Zilber lemma is a finite Morley rank gadget; \emph{the present contribution shows that it isn't}.
\item
Most sources focus on the ring \emph{generated} by the action instead of going to the centraliser; exceptions are \cite{NNonassociative, MPPrimitive}.
\item
Only \cite{NNonassociative} discusses rings and makes the connection with Schur's Lemma.
%All involve  (\cite{MPPrimitive} does introduce the centraliser rather than the ring generated.)
\item
\cite[487]{MPPrimitive} notices resemblances between various linearisation results:
\begin{quote}
{There appear to be no immediate implications between this and the results recorded here, though it looks similar to Theorem~1.2.}
\end{quote}
\emph{The present contribution elucidates the desired relations.}
\item
My own interest in the topic started when I read \cite{NNonassociative} while preparing \cite{DRegard}. This resulted in a very partial version of the present results, in finite Morley rank and using indecomposability. Then Wagner shared numerous ideas (to be continued elsewhere \cite{DWEndogenies}).
\end{itemize}

\subsection{Fields in soluble groups}\label{s:soluble}

To some extent, the Schur-Zilber lemma is the poor relation of the following theorem \cite[Corollary p.~175]{ZSome} (currently undergoing generalisation by Wagner):
\begin{quote}
connected, non-nilpotent, soluble groups of finite Morley rank interpret infinite fields.
\end{quote}

I believe the significance of the latter has been exaggerated for three reasons.
\begin{enumerate}
\item
In the local analysis of simple groups of finite Morley rank, different soluble subquotients may interpret non-isomorphic fields. Since there are strongly minimal structures interpreting \emph{different} infinite fields \cite{HStrongly}, any field structure could be a false lead. (For more on how experts approach the Algebraicity Conjecture on simple groups of finite Morley rank, and the influence of finite group theory instead of pure model theory, see both Cherlin's and Poizat's contributions in the same volume.)
\item
Fields obtained by this method can have `bad' properties, typically non-minimal multiplicative group \cite{BHMWBoese}.
\item
The corollary focused on abstract groups and distracted us from doing representation theory (see Borovik's remarkable contribution in the same volume).
\end{enumerate}

\section{The Theorem}\label{S:theorem}

\S~\ref{s:definitions} contains all necessary definitions. \S~\ref{s:explaining} justifies the structure of the statement. \S~\ref{s:optimality} discusses optimality, \S~\ref{s:corollaries} gives corollaries, and \S~\ref{s:indecomposable} considers the relation to `indecomposable generation'.

The general version of the skew-field lemma is a double-centraliser theorem, repeated below. Alternative names could have been `bimodule theorem' or `double-centraliser linearisation'.

\begin{theorem*}[Zilber's skew-field lemma]
Work in a finite-dimensional theory.
Let $V$ be a definable, connected, abelian group and $S, T \leq \DefEnd(V)$ be two invariant rings of definable endomorphisms such that:
\begin{itemize}
\item
$V$ is irreducible as an $S$-module (viz.~in the definable, connected category);
\item
$C(S) = T$ and $C(T) = S$, with centralisers taken in $\DefEnd(V)$;
\item
$S$ and $T$ are infinite;
\item
$S$ \emph{or} $T$ is unbounded.
\end{itemize}
Then there is a definable skew-field $\bK$ such that $V \in \Vect{\bK}_{<\aleph_0}$; moreover $S \simeq \End(V\colon \Vect{\bK})$ and $T \simeq \bK \Id_V$ are definable.
\end{theorem*}

It would be interesting to recast this kind of double-centraliser result in the abstract ring $S\otimes T$, with no reference to $V$.

\subsection{Definitions}\label{s:definitions}

\begin{itemize}
\item
Connected: with no definable proper subgroup of finite index. (Since the context does not provide a \textsc{dcc}, not all definable groups have a connected component.)
\item
Bounded: which does not grow larger when taking larger models. (The algebraist may fix a saturated model with inaccessible cardinality 
and argue there; bounded then means small. Also see \cite{HKSaturated}.)
\item
Type-definable: a bounded intersection of definable sets.
\item
Invariant: a bounded union of type-definable sets. (The name comes from the action of the Galois group of a `large' model. \S~\ref{s:explaining} gives reasons for considering the invariant category instead of the definable one.)
%Is the invariant class closed under bounded unions/intersections, and complement???
\item
Irreducible: no non-trivial proper submodule---a submodule being definable \emph{and connected}. (This is weaker than usual algebraic simplicity, which would also exclude finite submodules. Model theory will handle those in its own way.)
\item
Finite-dimensional: which bears a reasonable dimension on interpretable sets. \cite{WDimensional} would say \emph{fine, integer-valued, finite-dimensional}. The definition is as follows.
\end{itemize}

\begin{definition*}[{\cite{WDimensional}}]
A theory $T$ is {\upshape [}fine, integer-valued{\upshape ]} finite-dimensional if there is a dimension function $\dim$ from the collection of all interpretable sets in models of $T$ to $\bN\cup\{-\infty\}$, satisfying for a formula $\varphi(x, y)$ and interpretable sets $X$ and $Y$:
\begin{itemize}
\item
Invariance: If $a \equiv a'$ then $\dim(\varphi(x, a)) = \dim(\varphi(x, a'))$.
\item
Algebraicity: $X$ is finite nonempty iff $\dim(X) = 0$, and $\dim(\emptyset) = - \infty$.
\item
Union: $\dim(X \cup Y) = \max\{\dim(X), \dim(Y)\}$.
\item
Fibration: If $f\colon X \to Y$ is an interpretable map such that $\dim(f^{-1}(y)) \geq d$ for all $y \in Y$, then $\dim(X) \geq \dim(Y) + d$.
\end{itemize}
\end{definition*}

The dimension extends to type-definable, and then to invariant sets; of course one should no longer expect nice additivity properties.

Except for a key `field definability lemma' (\S~\ref{s:indecomposable}) we shall use little from \cite{WDimensional}.
There is an \textsc{acc} and a \textsc{dcc} on definable, \emph{connected} subgroups.

\subsection{Explaining the statement}\label{s:explaining}

Our statement deviates from traditional versions in several respects, and we make three cases for three notions.

\paragraph{Skew-fields rather than fields.}

Schur's lemma produces a skew-field, and so does Zilber's model-theoretic version.
\begin{itemize}
\item
This went first unnoticed since $\aleph_1$-categorical skew-fields are commutative. (Answering a question of Macintyre's, proved by Cherlin and Shelah---see note on \cite[p.~139]{BNGroups}---and independently by Zilber \cite{ZGruppi}.)
\item
It is easy to construct, in tame geometry, so-called `quaternionic representations', where the Schur field is the skew-field of quaternions.
\item
Also, the subring $\generated{A} \leq \End(V)$ generated by a \emph{commutative} group action can be smaller than its Schur skew-field $C_{\End(V)}(A)$: classical focus on the former (as in most sources) captures only partial geometric information.
\end{itemize}

So skew-fields are naturally unavoidable.
(There remains the question of which skew-fields can arise in a finite-dimensional theory. Skew-fields abound in number theory, but arguably number theory is far from tame. One can also doubt that the more exotic objects constructed in \cite{CSkew} will be finite-dimensional. The bold would conjecture that infinite skew-fields in finite-dimensional theories are: commutative and real closed, commutative and algebraically closed, or quaternionic over a commutative real closed field.)

\paragraph{Rings rather than groups.}

Let $V$ be an abelian group; then $\End(V)$ is a ring.
This accounts for studying representations of \emph{rings}.
\begin{itemize}
\item
If $G \leq \Aut(V)$ is a definable acting group, the subring of $\End(V)$ it generates need not be definable (see `invariance' below). This may have baffled pioneers in the topic.
\item
Rings were long neglected after Zilber's seminal \cite{ZGruppi} (a remarkable exception being \cite{NNonassociative}). Going to the enveloping ring however gives powerful results, inaccessible to group-theoretic reasoning: see Borovik's contribution in the present volume.
\end{itemize}

\paragraph{Invariance rather than definability.}

%We still have to account for invariance instead of definability.

Leaving definability may have stopped first investigators of the matter; it is however salutary.

\begin{itemize}
\item
If $G \leq \Aut(V)$ is a definable group, then the generated subring $\generated{G} \leq \End(V)$ is $\bigvee$-definable; this is closer to definability than invariance is.
However (see `skew-fields' above), $\generated{G}$ does not capture enough geometric information. The double-centraliser $C(C(G)) \geq \generated{G}$ is more adapted to Schur-style arguments.
\item
So let $R \leq \End(V)$ be a definable ring. Then Schur's covariance ring $C_{\DefEnd(V)}(R)$ need not be definable, but it is invariant. And if $R$ itself is invariant, $C_{\DefEnd(V)}(R)$ is too.
%TODO prove it.
%Say $R = \bigcup X_i$ where $X_i = \bigcap Y_{i, j}$, all operations being bounded. Then $C(R) = \bigcap C(X_i)$. But what then?
\end{itemize}
So model-theoretic invariance arises as naturally as centralisers do.

\subsection{Optimality}\label{s:optimality}

\begin{itemize}
\item
Both $S$ and $T$ must be infinite.

Let $\bK$ be a pure algebraically closed field of positive characteristic $p$ and $V = \bK_+$, which is definably minimal. Now $\DefEnd(V)$ consists of quasi-$p$-polynomials, viz.~of all maps $x \mapsto \sum_{k = -n}^n a_{p^k} \Fr_{p^k}$, where $\Fr$ is the Frobenius automorphism of relevant power, and $a_{p^k} \in \bK$; there is no bound on $n$. Only the action of $\bF_p$ commutes to all these.
We then let $S = \DefEnd(V)$ and $T = \bF_p$ (or vice-versa). The first is not definable.
\item
At least one must be unbounded.

Same $V$; now let $S$ be the ring of all quasi-$p$-polynomials \emph{with coefficients in $\bF_p$}, viz.~the subring of $\DefEnd(V)$ generated by $\Fr_p$ and its inverse. Then one easily sees that $C(S) = S$ is countable, and not definable.
\end{itemize}
% \end{remark}

\subsection{Corollaries}\label{s:corollaries}

I give three corollaries, proved in \S~\ref{s:proofcorollary}. The first relates the main, `double-centraliser' theorem to Schur's Lemma. The second retrieves what is called `Zilber's Field Theorem' in sources such as \cite{BNGroups}. The third is a variation coming from Nesin's work and isolated by Poizat.

\begin{corollary}[Schur-Zilber, one-sided form]\label{c:SZ}
Work in a finite-dimensional theory. Let $V$ be a definable, connected, abelian group and $S \leq \DefEnd(V)$ be an invariant, unbounded ring of definable endomorphisms. Suppose that $V$ is irreducible as an $S$-module.
Then $T = C_{\DefEnd(V)}(S)$ is a definable skew-field.
\end{corollary}
Corollary~\ref{c:SZ} is however not equivalent to our main result, which also covers the case of unbounded $T$ and infinite $S$.

\begin{corollary}[{see \cite[Théorème~IV.1]{DRegard}}]\label{c:D}
% A typical use of this theorem is: with same $V$, suppose $G \leq \DefAut(V)$ is a definable group such that $V$ is irreducible as a $G$-module \emph{and} $C_{\DefEnd(V)}(G)$ is infinite.
Work in a finite-dimensional theory. Let $V$ be a definable, connected, abelian group and 
$G \leq \DefAut(V)$ is a definable group such that $V$ is irreducible as a $G$-module \emph{and} $C_{\DefEnd(V)}(G)$ is infinite. Then $T = C_{\DefEnd(V)}(G)$ is a definable skew-field (so the action of $G$ is linear).
\end{corollary}

Corollary~\ref{c:D} (or a minor variation) unifies and should replace various results such as \cite[Lemma~2]{ZSome}, \cite[Theorem~4]{LWCanada}, \cite[Lemma~12]{NNonassociative}, \cite[Theorem~1.2.b]{MPPrimitive}, \cite[Théorème~IV.1]{DRegard}, \cite[Theorem~2.6]{PPSSimple}, \cite[Proposition~4.1]{MMTPermutation}. However there are no claims on finite generation.

% Another typical use is the following, proved in \S~\ref{s:proofcorollary}.

\begin{corollary}[after Nesin and Poizat]\label{c:NP}
Work in a finite-dimensional theory. Let $V$ be a definable, connected, abelian group and $R \leq \DefEnd(V)$ be an invariant, unbounded, commutative ring of definable endomorphisms. Suppose there is an invariant group $G \leq \DefAut(V)$ such that:
\begin{itemize}
\item
$V$ is irreducible as a $G$-module;
\item
$G$ normalises $R$;
\item
$G$ is connected.
\end{itemize}
Then there is a definable skew-field $\bK$ such that $V \in \Vect{\bK}_{<\aleph_0}$; moreover $R \hookrightarrow \bK\Id_V$ and $G \hookrightarrow \GL(V\colon \Vect{\bK})$. 
\end{corollary}

It would be interesting to relax the assumption on commutativity of $R$.
Further generalisations are expected using endogenies instead of endomorphisms \cite{DWEndogenies}.

\subsection{Indecomposable generation (and how to avoid it)}\label{s:indecomposable}

Contrary to widespread belief, the Schur-Zilber lemma has nothing to do with another celebrated result from Boris' early work: the `indecomposability theorem' \cite[Teorema 3.3]{ZGruppi}, which by analogy with the algebraic case I prefer to call the \emph{Chevalley-Zilber generation lemma} (with renewed hope that Boris will not mind being in good company). For more on the topic, see \S~8 of Poizat's contribution in the present volume.

Both results are often presented jointly, which serves neither clarity nor purity of methods.
In contrast, the proof given here relies on another phenomenon.

\begin{lemma*}[{field definability; extracted from \cite[Proposition~3.6]{WDimensional}}]
Work in a finite-dimensional theory. Let $\bK$ be an invariant skew-field such that:
\begin{itemize}
\item
there is an upper bound on dimensions of type-definable subsets of $\bK$;
\item
$\bK$ contains an invariant, unbounded subset.
\end{itemize}
Then $\bK$ is definable.
\end{lemma*}

The first clause is satisfied if there is a definable $\bK$-vector space of finite $\bK$-linear dimension.

\section{The Proofs}\label{S:proof}

% \subsection{Field Definability}
% 
% The proof relies on a key observation.
% 
% \begin{lemma*}[{field definability; extracted from \cite[Proposition~3.6]{WDimensional}}]
% Work in a dimensional theory. Let $\bK$ be a skew-field such that:
% \begin{itemize}
% \item
% there is an upper bound on dimensions of type-definable subsets of $\bK$;
% \item
% $\bK$ contains an invariant, unbounded subset.
% \end{itemize}
% Then $\bK$ is definable.
% \end{lemma*}
% \begin{proof}
% Let $X \subseteq \bK$ be type-definable with maximal dimension; by assumption $\dim X > 0$. Then for $a\in \bK$ consider the map :
% \[\begin{array}{ccc}
% X \times X & \to & \bK\\
% (x, y) & \mapsto & ax + y.
% \end{array}\]
% Since the image is type-definable and $\dim X^2 = 2 \dim X > \dim X$, the map is not injective. So there are $(x_1, y_1) \neq (x_2, y_2)$ with $ax_1 + y_1 = a x_2 + y_2$, or $a = \frac{y_2 - y_1}{x_2 - x_1}$. Hence $\bK$ is interpretable in $X^4$; now a type-definable field, it is definable.
% \end{proof}

The corollaries are derived in \S~\ref{s:proofcorollary}.
Let $V, S, T$ be as in the theorem. The proof is a series of claims arranged in Propositions.

\begin{proof}[Proof of Zilber's skew-field lemma]
It is convenient to let $T$ act from the right and treat $V$ as an $(S, T)$-bimodule.

\begin{proposition}\label{p:Tdomain}\leavevmode
\begin{enumerate}[label={\itshape (\roman*)},series=claims]
\item\label{i:Tdomain}
$T$ is a domain acting by surjections with finite kernels.
\end{enumerate}
\end{proposition}

This will later be reinforced in~\ref{i:Tskewfield}.

\begin{proof}\leavevmode
\begin{enumerate}[label={\itshape (\roman*)},series=proofs]
\item
Let $t \in T\setminus\{0\}$. Then $0 < Vt$ is $S$-invariant, definable, and connected; by $S$-irreducibility $Vt = V$ so $t$ is onto. In particular $T$ is a domain. Finally $\dim \ker t = \dim V - \dim Vt = 0$, so $\ker t$ is finite.
\qedhere
\end{enumerate}
\end{proof}

The global behaviour is difficult to control, so we go down to a more `local' scale with a suitable notion of lines.

\subsection{Lines}

\begin{notation}
Let $\delta = \min \{\dim sV: s \in S \setminus\{0\}\}$ and $\Lambda = \{s V: \dim sV = \delta\}$ be the set of \emph{lines}.
\end{notation}

\begin{proposition}\label{p:Vsum}\leavevmode
\begin{enumerate}[resume*=claims]
\item\label{i:LTinvariant}
Every line is $T$-invariant.
\item\label{i:sLinLambda}
If $L \in \Lambda$ and $s \in S$ are such that $sL \neq 0$, then $sL \in \Lambda$; in particular $L \cap \ker s$ is finite.
\item\label{i:Vsum}
$V$ is a finite sum of lines.
\item\label{i:Stransitive}
$S$ is transitive on $\Lambda$.
\end{enumerate}
\end{proposition}

Ads~\ref{i:sLinLambda} and~\ref{i:Vsum} will later be reinforced in~\ref{i:Lcapkers} and~\ref{i:Vdirectsum}, respectively.

\begin{proof}\leavevmode
\begin{enumerate}[resume*=proofs]
\item
Obvious since $S$ and $T$ commute.
\item
Say $L = s_0 V$. If $sL \neq 0$, then $0 < sL = (s s_0) V \leq s_0 V$, so by minimality of $\delta$ one has $sL \in \Lambda$. This also implies $\dim (L \cap \ker s) = \dim \ker s_{|L} = \dim L - \dim sL = 0$, and $L \cap \ker s$ is finite.
\item
The subgroup $0 < \sum \Lambda \leq V$ is definable, connected, and $S$-invariant; it equals $V$. Since dimension is finite, it is a finite sum.
\item
Let $L_1, L_2 \in \Lambda$, say $L_i = s_i V$. Now as above, $V = \sum_S s L_1 \not\leq \ker s_2$, so there is $s \in S$ such that $s_2 s L_1 \neq 0$. But then $0 < s_2 s L_1 = s_2 s s_1 V \leq s_2 V = L_2$, and equality holds.
\qedhere
\end{enumerate}
\end{proof}

\subsection{Linearising lines} 

\begin{proposition}\label{p:TonL}\leavevmode
\begin{enumerate}[resume*=claims]
\item\label{i:Lcapkers}
If $L \in \Lambda$ and $s \in S$ are such that $sL \neq 0$, then $L \cap \ker s = 0$.
\item\label{i:TonsomeL}
$T$ acts by automorphisms on \emph{some} line.
\item\label{i:ToneachL}
$T$ acts by automorphisms on \emph{each} line.
\end{enumerate}
\end{proposition}

The proof is different depending on whether $S$ or $T$ is unbounded.

\begin{proof}[Proof if $T$ is unbounded]\leavevmode
\begin{enumerate}[resume*=proofs]
\item
Suppose $sL \neq 0$; we show $L \cap \ker s = 0$. By~\ref{i:Stransitive}, $S$ is transitive on $\Lambda$, so there is $s' \in S$ with $s' s L = L$. Now $L \cap \ker s \leq L \cap \ker (s's)$, so we may assume that $s L = L$. Recall that $\ker s_{|L} = L \cap \ker s$ is finite by~\ref{i:sLinLambda}. Considering $s_{|L}^2\colon L \to L$ which is onto, we inductively find $|\ker s_{|L}^n| = |\ker s_{|L}|^n$, so $K = \sum_{n \in \bN} \ker s_{|L}^n$ is either trivial or countable infinite. \emph{Since $T$ is unbounded}, there is $t \in T \setminus\{0\}$ annihilating $K$. But $t$ has a a finite kernel by~\ref{i:Tdomain}, so $K = 0$, as desired.
\item
Let $t \in T\setminus\{0\}$ such that $\ker t > 0$. This finite subgroup is $S$-invariant and $S$ is infinite, so there is $s \in S\setminus\{0\}$ annihilating $\ker t$. Now $V = \sum \Lambda$ by~\ref{i:Vsum}, so there is $L_0 \in \Lambda$ with $sL_0 \neq 0$. Then $L_0 \cap \ker t \leq L_0 \cap \ker s = 0$ by~\ref{i:Lcapkers}.
\item
If $L_0$ is as in~\ref{i:TonsomeL} and $L$ is another line, by transitivity~\ref{i:Stransitive} there is $s \in S$ with $sL = L_0$. Then $s(L\cap \ker t) \leq L_0 \cap \ker t = 0$ so $L\cap \ker t \leq L \cap \ker s = 0$ by~\ref{i:Lcapkers}.
\qedhere
\end{enumerate}
\end{proof}

\begin{proof}[Proof if $S$ is unbounded]
The strategy is different here.
\begin{itemize}
\item
We first prove~\ref{i:TonsomeL}. By~\ref{i:Vsum} $V = \sum \Lambda$ is a finite sum, so there are $L_1, \dots, L_n$ such that $\bigcap_{i = 1}^n \Ann_S(L_i) = 0$. In particular $(S, +) \hookrightarrow S/\Ann_S(L_i)$ as abelian groups. \emph{Since $S$ is unbounded}, there exists some line $L$ such that the quotient group $\Sigma = S/\Ann_S(L)$ is unbounded. Let $t \in T\setminus\{0\}$. Then $K = \sum_{n \in \bN} \ker t_{|L}^n$ is either trivial or countable infinite. Since $\Sigma$ is unbounded, there is $\sigma \in \Sigma\setminus\{0\}$ annihilating $K$, i.e.~there is $s \in S$ annihilating $K$ but not $L$. By~\ref{i:sLinLambda} this shows $K = 0$, as desired.
\item
If $T$ acts on automorphisms on $L$, then for $s \in S$ with $s L \neq 0$ one has $L \cap \ker s = 0$. Indeed, $L \cap \ker s$ is finite by~\ref{i:sLinLambda}. Since $T$ is infinite there is $t \in T\setminus\{0\}$ with $(L \cap \ker s)t = 0$, but $t$ induces an automorphism of $L$. This proves~\ref{i:Lcapkers}, but only for lines on which $T$ acts by automorphisms.
\item
We deduce~\ref{i:ToneachL}, and~\ref{i:Lcapkers}.
Let $L$ on which $T$ acts by automorphisms and $L'$ be another line. Then by transitivity~\ref{i:Stransitive}, there is $s \in S$ with $sL = L'$. Suppose $w \in L' \cap \ker t$. Then there is $v \in L$ with $sv = w$. Now $s(vt) = (sv)t = wt = 0$, so $vt \in L \cap \ker s = 0$. Since $T$ acts by automorphisms on $L$, we have $v = 0$ and 
$w = 0$, as desired.
\qedhere
\end{itemize}
\end{proof}

Since it is unclear at this stage whether every element belongs to a line, we cannot immediately conclude that $T$ acts by automorphisms; this requires writing $V$ as a direct sum.

\subsection{Globalising local geometries}

\begin{proposition}\label{p:Tskewfield}\leavevmode
\begin{enumerate}[resume*=claims]
\item\label{i:linescomplemented}
Lines are complemented as $T$-modules, viz.~for $L \in \Lambda$ there is a definable, connected, $T$-invariant $H \leq V$ with $V = L \oplus H$.
\item\label{i:Vdirectsum}
$V$ is a finite, \emph{direct} sum of lines.
\item\label{i:Tskewfield}
$T$ is a skew-field acting by automorphisms.
\end{enumerate}
\end{proposition}

\begin{proof}\leavevmode
\begin{enumerate}[resume*=proofs]
\item 
% From~\ref{i:Lcapkers} we deduce~\ref{i:linescomplemented}. 
Say $L = s_0 V$. Since $V = \sum_S s L$ by~\ref{i:Vsum} and~\ref{i:Stransitive}, there is $s \in S$ with $s_0 s L \neq 0$, so $0 < s_0 s L = s_0 s s_0 V \leq L$. Let $s_1 = s_0 s$, so that $L = s_1 V = s_1 L$. Then $H = \ker s_1$ is such that $V = L + H$.
Now $L \cap H = L \cap \ker s_1 = 0$ by~\ref{i:Lcapkers}, so actually $V = L \oplus H$. Connectedness of $H$ follows.
\item
% In turn, \ref{i:linescomplemented} implies~\ref{i:Vdirectsum}.
Suppose that by~\ref{i:linescomplemented} lines $L_1, \dots, L_i$ are given with direct complements $H_j$ such that:
\begin{quote}
for $j \leq i$, one has $L_j \leq \bigcap_{k < j} H_k$ (viz.~each new line is in all previous complements).
\end{quote}
Then a quick induction yields:
\[V = \left(\bigoplus_{j = 1}^i L_j\right) \oplus \left(\bigcap_{k = 1}^i H_j\right).\]
Let $q$ project $V$ onto $\bigcap_{j = 1}^i H_j$ with kernel $\bigoplus_{j = 1}^i L_j$; it is $T$-covariant so $q \in C(T) = S$.

If $\bigoplus_{j = 1}^i L_j < V$, then $q \neq 0$; now $V = \sum \Lambda$ so there is $L' \in \Lambda$ such that $q L' \neq 0$. But then $L_{i+1} = qL' \in \Lambda$ and $L_{i+1} \leq \bigcap_{j = 1}^i H_j$ and we have iterated the process.
By finite-dimensionality it must terminate; and when it does, $V = \bigoplus_{i = 1}^n L_i$.
\item
Say $V = \bigoplus_{i = 1}^n L_i$ by~\ref{i:Vdirectsum}. Then for $t \in T$ one has $\ker t = \bigoplus_{i = 1}^n (L_i \cap \ker t) = 0$ by~\ref{i:ToneachL}.
\qedhere
\end{enumerate}
\end{proof}

Hence $T$ is a skew-field and $V \in \Vect{T}$, but we still fall short of definability. 

\subsection{Definability}

We return to lines. The next result is of a purely auxiliary nature.

\begin{proposition}\leavevmode
\begin{enumerate}[resume*=claims]
\item\label{i:localinverse}
Let $L_1, L_2 \in \Lambda$. If $\sigma \colon L_1 \simeq L_2$ is definable and $T$-covariant, then there is an \emph{invertible} $s \in S^\times$ inducing $\sigma$.
\end{enumerate}
\end{proposition}
\begin{proof}\leavevmode
\begin{enumerate}[resume*=proofs]
\item
Using~\ref{i:linescomplemented}, write $V = L_1 \oplus H_1$ for some $\pi_1 \in S$ with $L_1 = \im \pi_1$ and $H_1 = \ker \pi_1$.

If $L_2 \cap H_1 = 0$, then $H_1$ is a common direct complement for $L_1$ and $L_2$. Glue $\sigma\colon L_1 \to L_2$ with $\Id_H$ to produce a $T$-covariant map, viz.~an element of $C_{\DefEnd(V)}(T) = S$, inducing $\sigma$. It clearly is invertible.

If $L_2 \leq H_1$, then the process proving~\ref{i:Vdirectsum} enables us to take $L_1$ and $L_2$ as the first two lines in a direct sum decomposition. Consider the map given on $L_1$ by $\sigma$, on $L_2$ by $\sigma^{-1}$, and on the remaining sum by $1$. It is $T$-covariant and bijective, hence invertible in $S$; it induces $\sigma$.

The case $0 < L_2 \cap H_1 < L_2$ cannot happen. For then $\ker \pi_{1|L_2} \geq L_2 \cap H_1 > 0$ so by definition of lines, $\pi_1 L_2 = 0$ and $L_2 \leq H_1$.
\qedhere
\end{enumerate}
\end{proof}

\begin{notation}
For $L \in \Lambda$, by~\ref{i:linescomplemented} there exists a definable, connected, $T$-invariant $H$ such that $V = L \oplus H$.
\begin{itemize}
\item
Let $\pi_L$ be the relevant projection and $S_L = \pi_L S \pi_L$.
\item
Also let $T_L \leq \DefEnd(L)$ be the image of $T$.
\end{itemize}
\end{notation}

In full rigour, $S_L$ also depends on the complement chosen; we omit it from notation. This will not create difficulties.

\begin{proposition}\leavevmode
\begin{enumerate}[resume*=claims]
\item\label{i:SLTL}
$S_L$ and $T_L$ are skew-fields contained in $\DefEnd(L)$.
\item\label{i:CSLCTL}
Inside $\DefEnd(L)$ one has $C(S_L) = T_L$ and $C(T_L) = S_L$.
% \item\label{i:SLTL}
% $T \simeq T_L \simeq S_L^\op$ as rings.
\item\label{i:Tdefinable}
$T$ is definable.
\end{enumerate}
\end{proposition}
\begin{proof}
In case $T$ is unbounded, one may directly jump to~\ref{i:Tdefinable}.
\begin{enumerate}[resume*=proofs]
\item
Be careful that $S_L$ is an additive subgroup of $S$ closed under multiplication but it need not contain $1$. (Sometimes $S_L$ is called a \emph{subrng}, for `subring without identity'.)
However, $S_L$ \emph{per se} is a ring with identity $\pi_L$, and the latter acts on $L$ as $\Id_L$. Moreover if $\pi_L s \pi_L$ annihilates $L$, then since it annihilates the chosen direct complement, it is $0$ as an endomorphism of $V$, viz.~$\pi_L s \pi_L = 0$ in $S$. So $S_L$ can be viewed as a subring of $\DefEnd(L)$, and it is exactly the subring of restrictions-corestrictions $\{s_{|L}^{|L} : s \in \Stab_S(L)\}$. (This explains why the complement plays no role in our construction. It is however useful to have both points of view on $S_L$.)

Let $s \in S_L\setminus\{0\}$. Hence $sL = L$, so by~\ref{i:Lcapkers} it induces some $T$-covariant automorphism $\sigma$ of $L$; by~\ref{i:localinverse} there is $s' \in S$ inducing $\sigma$. Now $\pi_L s'^{-1} \pi_L$ is a two-sided inverse of $s$ in $S_L$.
This proves that $S_L$ is a skew-field; of course so is $T_L \simeq T$.
\item
One of them is easy. Let $f\colon L \to L$ be a definable, $T_L$-covariant morphism. By definition, $f$ commutes with the action of $T$. Take any $T$-invariant direct complement $H$ and set $\hat{f} = 0$ on $H$. Then $\hat{f}\colon V \to V$ is $T$-covariant. Hence $\hat{f} \in C(T) = S$ and $\pi_L \hat{f} \pi_L = f \in S_L$.

Now let $g\colon L \to L$ be definable and $S_L$-covariant. We aim at extending $g$ to an $S$-covariant endomorphism of $V$.

For $M \in \Lambda$ first choose $s \in S$ with $sL = M$ by transitivity~\ref{i:Stransitive}, then by~\ref{i:localinverse} we may assume $s \in S^\times$. Notice that $sgs^{-1}$ leaves $M$ invariant, and let $g_M\in \DefEnd(M)$ be the induced map.
We claim that this does not depend on the choice of $s$. Indeed let $s'$ be another invertible choice, giving rise to $g'_M$. Then $s^{-1} s'$ induces an element of $S_L$, so $g$ commutes with it and we find $g_M = g'_M$.

It also follows that $g_M \in C(S_M)$. For if $\eta \in S_M$ then we may assume $\eta \neq 0$ so by~\ref{i:localinverse} it is induced by an invertible element $h \in S^\times$ normalising $M$. Then $s' = hs$ is another invertible element taking $L$ to $M$. By the preceding paragraph, $s'gs'^{-1} = h g_M h^{-1}$ and $sgs^{-1} = g_M$ agree on $M$, so $g_M$ commutes with $\eta$.

We prove more: if $s\in S$ induces $\sigma\colon M \simeq N$, then $g_N \sigma = \sigma g_M$. Indeed by~\ref{i:localinverse}, freely supposing $s$ invertible, and picking invertible $s_M, s_N$ inducing $L \simeq M, N$, we find:
\[ g_N \sigma = s s^{-1} \cdot s_N g s_N^{-1} \cdot s \pi_M = s \cdot (s^{-1} s_N) g (s_N^{-1} s) \cdot \pi_M = s g_M.\]

Finally take a direct sum $V = \bigoplus L_i$ as in~\ref{i:Vdirectsum} and let $\hat{g} (\sum \ell_i) = \sum g_{L_i}(\ell_i)$, which is definable, well-defined, and extends $g$. We want to show $\hat{g} \in C(S)$. Let $s \in S$; also let $s_i = \pi_i s$. It is enough to show that $\hat{g}$ commutes with each $s_i$, and it is enough to show that they commute on each $L_j$.
We have thus reduced to checking that $\hat{g}$ and $\sigma \colon L_j \simeq L_i$ induced by an element of $S$ commute. But this is the previous paragraph.

Hence $\hat{g} \in C(S) = T$ and therefore $g = \hat{g}_{|L} \in T_L$.
\item
If $T$ is unbounded we directly apply the field definability lemma from \S~\ref{s:indecomposable} (in that case, \ref{i:SLTL} and~\ref{i:CSLCTL} are not necessary). So we suppose that $S$ is unbounded.

We first prove that there is $L$ such that $S_L$ is unbounded.
By~\ref{i:Vdirectsum} take any decomposition $V = \bigoplus_{i = 1}^n L_i$ and form projections $\pi_i$ onto $L_i$ with kernels $\bigoplus_{j \neq i} L_j$. Let $S_{i, j} = \pi_i S \pi_j$, an additive subgroup of $S$. We contend that one of them is unbounded. Indeed, the additive group homomorphism:
\[\begin{array}{ccc}
S & \to & \prod_{i, j} S_{i, j}\\
s & \mapsto & (\pi_i s \pi_j)_{i, j}
\end{array}\]
is injective since $\sum_k \pi_k = 1$.
Now if $S_{L, M}$ and $S_{L', M'}$ are defined as the $S_{i, j}$, one easily sees $S_{L, M} \simeq S_{L', M'}$ definably; so all rings $S_L$ are unbounded.

Caveat: because $S_L$ and $T_L$ are mutual centralisers only in $\DefEnd(L)$ and not in $\End(L)$, the following paragraph cannot be made more trivial.

Therefore $S_L$ is an unbounded skew-field by~\ref{i:SLTL}. By field definability of \S~\ref{s:indecomposable}, $S_L$ is definable; now $\dim S_L > 0$ and $\dim L$ is finite, so $L\in \Vect{S_L}_{< \aleph_0}$. In particular, all $S_L$-endomorphisms of $L$ are definable, so by~\ref{i:CSLCTL} one has $T_L = \End(L\colon \Vect{S_L})$. This is a skew-field by~\ref{i:SLTL}, so the \emph{linear} dimension over $S_L$ is $1$ and $T\simeq T_L \simeq S_L^\op$ is unbounded as well.
\qedhere
\end{enumerate}
\end{proof}

By field definability, the skew-field $T$ is definable infinite, so $\dim T > 0$; now $\dim V$ is finite so $V \in \Vect{T}_{< \aleph_0}$. Finally $S = C(T) = \End(V\colon\Vect{T})$. Lines in our sense now coincide with $1$-dimensional $T$-subspaces of $V$.
This completes the proof of Zilber's skew-field lemma.
\end{proof}

\subsection{Proofs of Corollaries}\label{s:proofcorollary}

\begin{proof}[Proof of Corollary~\ref{c:SZ}]
Replacing $S$ by $C_{\DefEnd(V)}(T)$ preserves irreducibility, invariance, unboundedness, and $T$. Notice that the latter invariant ring acts by surjective endomorphisms, so is a domain. If $T$ is a finite domain, then it is a definable field. If $T$ is infinite, apply the Theorem.
\end{proof}

\begin{proof}[Proof of Corollary~\ref{c:D}]
Let $T = C_{\DefEnd(V)}(G)$ and $S = C_{\DefEnd(V)}(T) \supseteq G$. Apply the Theorem.
\end{proof}

\begin{proof}[Proof of Corollary~\ref{c:NP}]
Let $V, R, G$ be as in the statement.
The proof follows that of \cite[Théorème~3.8]{PGroupes} closely.
% Let $S = C_{\DefEnd(V)}(R)$, which contains $R$ by commutativity.
Let $W \leq V$ be $R$-irreducible, viz.~minimal as a definable, connected, $R$-submodule; this exists by the \textsc{dcc} on definable, connected subgroups. Let $\fp = \Ann_R(W)$, a relatively definable ideal of $R$.

For $g \in G$, the definable, connected subgroup $g W \leq V$ is $R$-invariant: hence an $R$-submodule. Clearly $\Ann_R(gW) = g \fp g^{-1}$. Moreover $R/\fp \simeq R/(g\fp g^{-1})$.

Now by $G$-irreducibility, $V = \sum_G gW$. So there are $g_1, \dots, g_n \in G$ such that $V = \sum_{i = 1}^n g_i W$. In particular $\bigcap_{i = 1}^n \Ann_R(g_i W) = 0$, and $R \hookrightarrow \prod R/(g_i \fp g_i^{-1})$. We just saw that all terms have the same cardinality. They are therefore unbounded.

Hence the unbounded, commutative ring $R/\fp$ acts faithfully on the $R/\fp$-irreducible module $W$. Notice that $R/\fp \leq C_{\DefEnd(W)}(R/\fp)$. By the Theorem, the action of $R/\fp$ on $W$ is linearisable, and $R/\fp$ acts by scalars.
The problem is to make this linear structure global without losing the action of $G$.
But we know that $\fp$ is a prime ideal of $R$.

Now consider the set of prime ideals $P = \{h\fp h^{-1}: h \in G\}$. Suppose $\fp_1, \dots, \fp_k \in P$ are distinct, say $\fp_i = h_i \fp h_i^{-1}$. By prime avoidance, there are elements $r_i \in \fp_i \setminus\bigcup_{j \neq i} \fp_j$. Then taking products, there are elements $r'_i \in \bigcap_{j \neq i} \fp_j \setminus \fp_i$. These are used to show that the sum $\sum_{i = 1}^k h_i W$ is direct. In particular, $k \leq \dim V$ and $P$ is finite.

Since $G$ is connected and transitive on the finite set $P$, the latter is a singleton, namely $P = \{\fp\}$. But by faithfulness one had $\bigcap P = 0$, so $\fp = 0$.

Now let $r \in R\setminus\{0\}$. Then $r \notin \fp$ acts on $W$ as a non-zero scalar, so $W \leq \im r$. Since $r$ was arbitrary, for any $g \in G$, one has $gW \leq \im r$. Summing, $\im r = V$; this implies that $\ker r$ is finite. Then $K = \sum_{n \in \bN} \ker r^n$ is either trivial or countable infinite. But by commutativity, it is $R$-invariant. Since $R$ is unbounded, there is $r_0 \in R\setminus\{0\}$ annihilating $K$. Since $r_0$ has a finite kernel in $V$, we see $K = 0$.
Thus the domain $R$ acts by automorphisms on $V$.

Hence $\bF = \Frac(R)$ is naturally a subring of $\DefEnd(V)$. By field definability, it is definable. Now $G$ normalises $\bF$ and centralises it \cite[\S~3.3]{WDimensional}. In particular, $G$ centralises $R$. Therefore  
$S = C_{\DefEnd(V)}(R)$, which contains $R$ by commutativity, also contains $G$. It follows that $V$ is $S$-irreducible and we apply the Theorem globally to conclude.
\end{proof}

% \section*{Acknowledgements}
\paragraph*{Acknowledgements.}
It is a pleasure to thank, alphabetically:
Tuna Alt{\i}nel, who first taught me Zilber's field theorem; Alexandre Borovik and Gregory Cherlin, who took me further into model-theoretic algebra; Ali Nesin and Bruno Poizat, for sharing their culture and passion; Frank Wagner, to whom I have a great technical debt; and last but not least, for a long-lasting vision which inspired so many of us, Boris Zilber.
%Boris Zilber for his inspiring vision, Frank Wagner to whom I have a great technical debt, and for their historical insight: Simon Thomas, Bruno Poizat, Ali Nesin, Gregory Cherlin and Alexandre Borovik.

\printbibliography

\end{document}